\documentclass[12pt]{amsart}

\usepackage[utf8]{inputenc}
\usepackage{amstext,amsthm,amsmath}

\usepackage{latexsym}
\usepackage{amsfonts}
\usepackage{mathtools}
\usepackage{graphicx}
\usepackage[normalem]{ulem}
\usepackage{fullpage}
\usepackage{marvosym}
\usepackage[usenames,dvipsnames]{color}
\usepackage{verbatim}
\usepackage[hidelinks]{hyperref}
\usepackage[english]{babel}
\usepackage{tabularx}
\usepackage{tikz}
\usepackage{enumerate}

\usepackage{bigints}

\newtheorem*{theorem*}{Theorem}
\newtheorem*{lemma*}{Lemma}
\newtheorem*{question*}{Question}
\newtheorem*{conjecture*}{Conjecture}
\newtheorem*{open problem*}{Open Problem}
\newtheorem*{problem*}{Problem}
\newtheorem*{proposition*}{Remark}
\newtheorem{theorem}{Theorem}

\newtheorem{lemma}[theorem]{Lemma}

\theoremstyle{definition}\newtheorem{remark}[theorem]{Remark}
\theoremstyle{definition}

\allowdisplaybreaks

\renewcommand\d{\Delta}
\newcommand\Hol{{\mathrm{Hol}}}

\newcommand\CC{{\mathbb C}}
\newcommand\RR{{\mathbb R}}

\newcommand\BB{{\mathbb B_2}}
\newcommand\SSS{{\mathbb S_2}}

\newcommand\Ba{\mathrm{B}}

\title[boundary zeroes]{Characterization of the boundary zero set of stable polynomials in $\mathbb B_2$}
\author[D. Vavitsas]{Dimitrios Vavitsas}
\thanks{}
\address{School of Mathematics (Zhuhai), Sun Yat-Sen University, Zhuhai, Guangdong, 519082, P. R.
China}
\email{vavitsas@mail.sysu.edu.cn}


\begin{document}

\begin{abstract}
Given a polynomial $p\in \mathbb{C}[z,w]$ whose zeros lie outside the unit ball $\mathbb B_2$ of $\CC^2,$ we study its local Weierstrass polynomial--obtained via the Weierstrass Preparation Theorem--and provide a strong version of the Nash curve selection lemma. These results lead to a characterization of the intersection of the zero set of $p$ with the unit sphere: $\mathcal{Z}(p)\cap\mathbb{S}_{2}.$
\end{abstract}

\maketitle

\section{Introduction}\label{Intro}

This paper is about the study of the zero set of polynomials in $\CC^2$ that satisfy a stability condition with respect to the unit ball.

In particular, we obtain a stronger version of the so called Nash curve selection lemma for polynomials $p\in \CC[z,w]$ whose zeros lie outside the Euclidean unit ball $\BB$ of $\CC^2$ and  achieve a characterization of the intersection of the zero set of $p$ with the unit sphere: $\mathcal{Z}(p)\cap\mathbb{S}_{2}.$ Such polynomials are often called stable polynomials, and the set of their zeros lying on the unit sphere is called boundary zero set.

These results also serve as a different  approach to recent breakthroughs by G. Knese, J.E. Pascoe, and A. Sola \cite[Proposition 1.3, Theorem 1.5]{Knese}, who gave a complete local characterization and parameterization of the (complex) algebraic curves—and their finitely many Newton-Puiseux branches—contained within the local zero set of a stable polynomial in $\CC^2.$

The \emph{Nash curve selection lemma} is a foundational result in semi-algebraic geometry; consequently, establishing stronger variants is a matter of considerable importance. Recall the injective version as it is stated by K. Kurdyka in \cite{Kurdyka1}. See also \cite[Proposition 8.1.13, Theorem 9.7.10]{Algebraic} which is a semi-algebraic consequence of the statements formulated in \cite{Cartan} by F. Bruhat, H. Cartan and \cite{Milnor} by J. Milnor. 

\begin{lemma*}
Let $A\subset \RR^N$ be a semi-algebraic set and let $\alpha\in \overline{A}$ be an accumulation point of $A.$ Then there exists an injective Nash curve $\gamma_\alpha:(-1,1)\rightarrow\RR^N$ such that 
$$\gamma_\alpha((0,1))\subset A \text{ and } \gamma_\alpha(0)=\alpha.$$
\end{lemma*}

It is possible that for each such curve $\gamma_\alpha$ and each $\epsilon>0,$ the set $\gamma_\alpha((-\epsilon,\epsilon))$ is not a (smooth) submanifold in $\RR^N;$ for instance consider the algebraic set $A:=\{(t^2,t^3):t\in \RR\}$ in $\RR^2$ and the point $\alpha=(0,0).$

Given a polynomial $p\in \CC[z,w],$ the set $\mathcal{Z}(p)\cap \SSS$ may be regarded as an algebraic set in $\RR^4.$ Hence, it admits a decomposition $\mathcal{Z}(p)\cap \SSS=\cup_iM_i,$ where $M_i$ are finitely many, pairwise disjoint Nash submanifolds of real dimension at most 1, see \cite[Theorem 1]{VZ}. In particular, the decomposition $\cup_iM_i$ may be chosen so that each one-dimensional component $M_i$ has exactly two points on its boundary, see the discussion right after \cite[Proposition 4.1]{VWZ}. According to \cite{Knese} and \cite{VZ}, we are led to the following question:

\begin{question*}
Given a stable polynomial $p\in \CC[z,w]$ with $\mathcal{Z}(p)\cap \SSS=\cup_iM_i,$ and any fixed accumulation point $\alpha\in \overline{M_i}\setminus M_i,$ does there always exist a Nash curve $\omega_\alpha:(-1,1)\rightarrow\mathcal{Z}(p)\cap \SSS,$ that is a diffeomorphism onto its image and satisfies $\omega_\alpha((0,1))\subset M_i,$ $\omega_\alpha(0)=\alpha$? 
\end{question*}

We shall prove that the answer to the  \emph{Question} is \emph{yes} thereby characterizing the boundary zero set $\mathcal{Z}(p)\cap \SSS.$

The main goal is to apply Weierstrass preparation theorem and Puiseux Parametrization theorem, together with arguments and ideas coming up in \cite{HW} by F.R. Harvey and R.O. Wells, \cite{Fornaess} by J.E. Fornaess and M. Ovrelid, \cite[Lemma C.3]{KneseLondon} by G. Knese, \cite{KneseK} by G. Knese  \L. Kosi\'{n}ski,  T. Ransford and A. Sola, \cite{Knese} by G. Knese, J.E. Pascoe, A. Sola, \cite{KV} by \L. Kosi\'{n}ski et al., \cite{Noel} by A. Noell, B. Stensønes, \cite[Proposition 10.7.5-Theorem 10.7.6]{Rudin} by W. Rudin, \cite{VWZ} by J. Wu, K. Zarvalis et al., and \cite{VZ} by K. Zarvalis et al.. 

The answer to the \emph{Question} is \emph{no} in the setting $\mathbb B_n,$ for $n\geq 3$ (see \cite{VWZ}); hence, the ambient $\CC^2$ is essential. 

The aforementioned works  provide the primary motivation for the present study. Readers may certainly consider natural extensions of this problem such as changing the domain of stability to one where the topological boundary is an algebraic set or real analytic set.
 
\section{Preliminary Remarks }\label{Prel Remarks}

Denote the \emph{unit ball} of $\CC^2$ by $\BB:=\{(z,w)\in \CC^2:||(z,w)||<1\}$ and its topological boundary, the \emph{unit sphere}, by $\SSS:=\{(z,w)\in \CC^2:||(z,w)||=1\},$ where $||(z,w)||=\sqrt{|z|^2+|w|^2}$ is the \textit{Euclidean norm} induced by the \textit{Euclidean inner product} $\langle (z,w),(u,v)\rangle=z\overline{u}+w\overline{v}$ in $\CC^2.$ Denote the disk with center $z\in \CC$ and radius $\epsilon>0$ by $\d(z,\epsilon):=\{z'\in\CC:|z-z'|<\epsilon\}.$

For the zero set of a function $f\in \CC[z,w],$ we use the notation $\mathcal{Z}(f)\vcentcolon=\{(z,w)\in \CC^2:f(z,w)=0\}$. In particular, a polynomial $p\in \CC[z,w]$ is said to be \emph{stable} (with respect to the unit ball) if $p$ is zero-free in the unit ball, that is, $\mathcal{Z}(p)\cap \BB=\emptyset$.  The set of the zeros of a stable polynomial that lie on the unit sphere is called \emph{boundary zero set}, and it  is written as  $\mathcal{Z}(p)\cap \SSS$.

Any set $A\subset \CC^n$ may be regarded as a subset of $\RR^{2n},$ thus we define its \emph{real dimension} by
$$\textup{dim}_\RR A\vcentcolon=\max\{\textup{dim}(\Gamma):\Gamma\subset A, \Gamma \text{ submanifold of }\RR^{2n}\}\in\{-1,0,1,...,2n\}$$
The case $\textup{dim}_\RR A=-1$ is devoted to the instance when $A=\emptyset$.

Recall the following definitions on semi-algebraic geometry; see \cite{Algebraic} for a concise presentation. See also \cite{Krantz} on theory of real analytic functions. Fix $N\in \mathbb N.$ A set $A\subset\RR^N$ is said to be \textit{semi-algebraic set} provided there exists a finite number of real polynomials $p_i,$ $i=1,...,I,$ and $q_{i,j}$,  $j=1,\dots,J$ such that $$A=\cup_{i=1}^{I}\cap_{j=1}^{J}\{x\in \RR^N:p_i(x)=0,q_{i,j}(x)>0\}.$$

\begin{itemize}
\item[(i)] Let $A\subset\RR^M$ and $B\subset\RR^N$ be two semi-algebraic sets. A mapping $f\vcentcolon A\to B$ is characterized as \textit{semi-algebraic} if its graph is a semi-algebraic set of $\RR^{M+N}$.
\item[(ii)] Let $A\subset\RR^N$ be semi-algebraic. A semi-algebraic mapping $f\vcentcolon A\to\RR$ of class $C^\infty$ is called a \textit{Nash function}. Moreover, given two semi-algebraic sets $A,B\subset\RR^N$, a semi-algebraic bijection $f\vcentcolon A\to B$ of class $C^\infty$ is called a \textit{Nash diffeomorphism}.
\item[(iii)] A semi-algebraic subset $M$ of $\RR^N$ is said to be a \textit{Nash submanifold} of $\RR^N$ of dimension $d$ if for every point $x\in M$, there exists a Nash diffeomorphism $\phi$ from an open semi-algebraic neighborhood $\Omega$ of the origin in $\RR^N$ onto an open semi-algebraic neighborhood $\Omega'$ of $x$ in $\RR^N$ such that $\phi(0)=x$ and $\phi((\RR^d\times\{0\})\cap\Omega)=M\cap\Omega'$.
\item[(iv)] Two Nash submanifolds $M$ and $N$ are considered to be \textit{Nash diffeomorphic} if there exists a bijection $f\vcentcolon M\to N$ such that both $f$ and $f^{-1}$ are Nash functions. The class of Nash functions coincide with the class of semi-algebraic real analytic functions.
\end{itemize}

Writing $p=\textrm{Re }p+i\textrm{Im }p,$ we see that $\mathcal{Z}(p)\cap \SSS$ is an semi-algebraic set in $\RR^4;$ indeed, it is the following intersection: $\{\textrm{Re }p=0\}\cap \{\textrm{Im }p=0\}\cap \{|z|^2+|w|^2=1\}.$ 

Recall the following definitions in the setting of the unit ball $\BB.$ We refer the reader to \cite{Rudin} for background on interpolation theory. Let $I\subset\RR$ be an interval on the real line. A $C^1$ curve $\gamma:I\to\SSS$ is said to be a \textit{complex tangential curve} when
$$\langle\gamma'(t),\gamma(t)\rangle=0, \quad\quad\textup{for all }t\in I.$$

According to \cite[Theorem 1]{VZ}, if $p\in \CC[z,w]$ is stable and $\dim_\RR(\mathcal{Z}(p)\cap \SSS)>0$, then $\mathcal{Z}(p)\cap \SSS=\cup_iM_i,$ where $M_i$ are pairwise disjoint Nash submanifolds, each equivalent to $(-1,1)^{m_i},$ $m_i=0,1,$ by a complex tangential Nash diffeomorphism. The term "equivalent" means that there exists a complex tangential Nash diffeomorphism $\phi_i:(-1,1)^{m_i}\rightarrow M_i$. By the discussion right after \cite[Proposition 4.1]{VWZ}, we may assume that each one-dimensional component $M_i$ contains exactly two points on its boundary $\overline{M_i}\setminus M_i$. Henceforth, the family $\{M_i\}_i$ refers to a such decomposition. If $\textrm{dim}_\RR(\mathcal{Z}(p)\cap \SSS)=0,$ then $\mathcal{Z}(p)\cap \SSS$ consists of finitely many points in the unit sphere.  Consequently, the more involved scenario arises when the boundary zero set contains at least one accumulation point. 

For the moment, the description provided above doesn't rule out the geometric case where the union of the closures of two distinct one-dimensional components, namely, $\overline{M_i}\cup\overline{M_j},$ $i\neq j,$ build a "cusp" at a accumulation point $\alpha\in\overline{M_i}\cap\overline{M_j}$.  At this point and throughout the paper, by the term "cusp" at $\alpha$, we mean that for $\alpha\in \overline{M_i}\setminus M_i,$ and each injective Nash curve $\gamma_\alpha:(-1,1)\rightarrow \RR^4$ derived from Nash curve selection lemma (see  Section~\ref{Intro}), there is no one-dimensional (smooth) submanifold $N\subset \RR^4$ such that $\gamma_\alpha([0,\epsilon))\subset N,$ for all $\epsilon>0.$

Let us suppose that $\dim_\RR(\mathcal{Z}(p)\cap \SSS)=1$ with $\mathcal{Z}(p)\cap \SSS=\cup_iM_i.$ In this specific case,  according to \cite{VWZ}, for each accumulation point $\alpha\in \overline{M_i}\setminus M_i$--for simplicity assume that $\alpha=(0,1)$; otherwise make use of a unitary matrix--there exists a Nash curve  $\gamma_\alpha:(-1,1)\rightarrow\mathcal{Z}(p)\cap \SSS,$ which is homeomorphism onto its image, complex tangential, and is classified by
$$\gamma_\alpha(t)=(\beta_\nu t^\nu+o(t^\nu),1-\frac{|\beta_\nu|^2}{2}t^{2\nu}+o(t^{2\nu})),$$ 
for some fixed integer $\nu\geq2$ and $\beta_\nu\neq 0.$ See \cite{Rudin} on unitary matrices. In particular, $\gamma_\alpha(0)=\alpha$ and $\gamma_\alpha'(0)=0;$ hence, removing the singularity at $\alpha$ via parameterization may not be possible. We shall show that this description can be further refined.  

For instance, consider the well known stable polynomial $p(z)=1-z^2-w^2$ and $\alpha=(0,1)\in \mathcal{Z}(p)\cap \SSS;$ the Nash curve $\gamma_{\alpha}(t)=(t^\nu,\sqrt{1-t^{2\nu}}),$ $t\in (-1,1),$ $\nu=3,5,...,$ is a homeomorphism onto its image, complex tangential and satisfies $\gamma_{\alpha}((-1,1))\subset \mathcal{Z}(p)\cap \SSS,$ $\gamma_\alpha(0)=\alpha.$ In particular, it can be parametrized as follows: $\gamma_\alpha(t)=(t,\sqrt{1-t^2})$ which yields a submanifold lying in $\mathcal{Z}(p)\cap \SSS.$. See \cite{Knese} for crucial and novel examples of stable polynomials.  

The results below examine whether the Nash curve $\gamma_\alpha$ may be chosen to satisfy $\gamma_\alpha'(0)\neq 0$.

\section{Main results}

The following lemma is a stronger version of the Nash curve selection lemma for the set $\mathcal{Z}(p)\cap \SSS$ of a stable polynomial $p\in \CC[z,w]$. The proof describes also the Weierstrass polynomial of the local Weierstrass factorization of the polynomial.

\begin{lemma}\label{Lemma}
Let $p\in \CC[z,w]$ be a stable polynomial with $\mathcal{Z}(p)\cap \SSS=\cup_iM_i,$ and let $\alpha\in \overline{M_i}\setminus M_i,$ for any fixed $i.$ There exists a Nash curve $\omega_\alpha:(-1,1)\rightarrow \mathcal{Z}(p)\cap \mathbb S_2$ that is a diffeomorphism onto its image and satisfies $\omega_\alpha((0,1))\subset M_i$, $\omega_\alpha(0)=\alpha.$

\begin{proof}
For simplicity assume that $p$ is irreducible and $\alpha:=(0,1).$ By the Nash curve selection lemma (see the Lemma in Section~\ref{Intro}), the identity principle and the arguments in the proof of \cite[Theorem 1]{VZ}, there exists a Nash curve  $\gamma_\alpha=(\gamma_1,\gamma_2):(-1,1)\rightarrow\mathcal{Z}(p)\cap \SSS,$ which is homeomorphism onto its image, complex tangential, and $\gamma_\alpha((0,1))\subset M_i,$ $\gamma_\alpha(0)=\alpha$.

Recall the discussion concerning the Weierstrass Preparation Theorem (see \cite[ 14.1.3. Theorem]{Rudin}) and the application of Puiseux's Parametrization Theorem (see \cite[pp. 136--137]{Fisher} and \cite[Corollary, pp. 171--171]{Loja}) in \cite[pp. 33-34]{KneseK} and \cite[349-350]{KV}. We may expand
$$p(z,w)=A(z,w)(1-h_1(z)w)\cdots (1-h_N(z)w), \quad (z,w)\in \d(0,\epsilon)\times \d(1,\epsilon),$$
where $A$ is holomorphic and zero free on the polydisk  $\d(0,\epsilon)\times \d(1,\epsilon),$ and each function $h_j$ is well defined, singled-valued holomorphic function on cut-branch subdomains of $\d(0,\epsilon)$ with respect to $0.$

In particular, there exists $\delta>0$, a non-negative integer $k>0,$ and holomorphic functions $\Phi_1,...,\Phi_N\in \Hol (\d(0,\delta))$ such that $h_j(z)=\Phi_{j(i)}(z^{1/k}),$ for a properly chosen branch and $z$ close to $0.$ Fix $j$ and $i(j):=i.$ Since $p$ is a stable polynomial, we obtain that 
\begin{equation}\label{Phi ineq stability 1}
|\Phi_i(z)|^2\leq \frac{1}{1-|z|^{2k}}, \quad \text{ or,} \quad |z|^{2k}+\frac{1}{|\Phi_i(z)|^2}-1\geq 0,
\end{equation}
for all $z$ close to $0.$ 

Then there exists $j$ (for simplicity, put $j=i(j)=i$) and $\epsilon>0$ such that
\begin{equation}\label{h Nash curve}
h_i(\gamma_1(t))=\frac{1}{\gamma_2(t)},
\end{equation}
for all $t\in (0,\epsilon)$ and fixed $\epsilon>0.$

The conditions above yield the following classification:
\begin{equation}\label{class h and Phi}
h_i(z)=\Phi_i(z^{1/k})=1+\lambda^{i}_{0} z^2+\lambda^i_{1}z^{2+1/k}+\lambda^{i}_{2}z^{2+2/k}....=1+\lambda^i_{0}z^2+o(z^2),
\end{equation}
where $|\lambda^i_0|=1/2,$ see \cite{KV}. For the moment, $h_i$ may be not holomorphic.

Let us suppose first that $h:=h_i$ is a holomorphic function on $\d(0,\epsilon),$ for some $\epsilon>0$ such that $1/h$ is well defined and holomorphic function on $\d(0,\epsilon)$ as well. This case can be handled utilizing the techniques in the proof of \cite[Proposition 4.6]{VWZ}; however, for completeness, we provide the proof. Indeed, define $\varphi:\d(0,\epsilon)\rightarrow\RR$ by
\begin{equation*}
    \varphi(z):=|z|^2+|h(z)|^{-2}-1
\end{equation*}
Stability of the polynomial yields $\varphi(z)\geq 0,$ for all $z.$ A straightforward computation shows that $\partial_{z\bar{z}}\varphi(0)>0.$ Whence, there exists a real analytic curve $\tilde{\omega}:(-1,1)\rightarrow \CC$ which is diffeomorphism onto its image, $\tilde{\omega}(0)=0,$ and there exists $\epsilon'>0$ such that $\mathcal{Z}(\varphi)\cap\d(0,\epsilon')\subset \tilde{\omega}((-1,1)),$ see \cite{HW} (or, the proof of \cite[Theorem 10.7.6]{Rudin}; apply the real analytic version of implicit function theorem, see \cite{Krantz}). Define $\omega_\alpha:(-1,1)\rightarrow \CC^2$ by $\omega_\alpha(t):=(\tilde{\omega}(t),1/h(\tilde{\omega}(t)).$ Then $\omega_\alpha$ is a real analytic curve, diffeomorphism onto its image, and $\omega_\alpha(0)=\alpha.$ Equation~\eqref{h Nash curve}  yields that $ \mathcal{Z}(\varphi)\cap \d(0,\epsilon)$ contains a specific set with an accumulation point, hence, $\tilde{\omega}((-1,1))$ contains this set too.  Summing up, $p(\omega_\alpha(t))=0,$ and $||\omega_\alpha(t)||^2=1,$ on a set having an accumulation point. Identity principle ensures that $\omega_\alpha:(-1,1)\rightarrow\mathcal{Z}(p)\cap \SSS.$ By the nature of the curves $\gamma_\alpha,$ $\omega_\alpha$, and the conditions that satisfy, we have $\omega_\alpha((0,\epsilon))\subset M_i$ for some fixed $\epsilon>0,$ and $\omega_\alpha(0)=\alpha.$ In particular, $\omega_\alpha((0,\epsilon))\cap M_i$ is connected and open in the relative topology of $M_i,$ $\alpha\in \overline{M_i\cap \omega_\alpha((0,\epsilon))},$ and $\omega_\alpha((-\epsilon,\epsilon))$ has exactly two points on its boundary.   Ergo, utilizing $\omega_\alpha$ we derive the desired Nash curve. 

We shall prove that stability forces $h$ to be always a holomorphic function and the assertion will follows from the previous step. The idea is to apply again the implicit function theorem (as in the proof of \cite[Theorem 10.7.6]{Rudin}) together with the stability condition in order to eliminate the terms that makes $h$ a non-holomorphic function. We slightly switch the stability inequality~\eqref{Phi ineq stability 1} to an equivalent one. 

Put $\Phi:=\Phi_i$ and $\lambda^i_{l}=\lambda_{l},$ $l=0,1,....$ If $k=1,$ then indeed it is holomorphic. Let us suppose in contrary that $k>1,$ and there exists at least one $l\notin\{0,k,2k,3k,...\}$ such that $\lambda_l\neq 0.$ Note that $\Phi(0)=1,$ hence, there exists a function $\Psi$ which is holomorphic close to $0,$ and such that $\Phi(z)=e^{\Psi(z)}.$ The inequality \eqref{Phi ineq stability 1} yields the equivalent one
\begin{equation}\label{real psi ineq sigma}
    \textrm{Re } \Psi(z)\leq \sigma(|z|^{2k}),
\end{equation}
for all $z$ close to $0,$ where $\sigma(\theta):= -\frac{1}{2}\log(1-\theta)=\frac{\theta}{2}+\frac{\theta^2}{4}+...,$ $|\theta|<1.$ In particular, $\Psi$ has a Taylor series expansion of the form
\begin{equation}\label{eq: Psi 1 classif}
\Psi(z)=\log(1-(1-\Phi(z))=- \sum_{q=1}^{+\infty}\frac{1}{q}(1-\Phi(z))^q
=\lambda_0 z^{2k}+o(z^{2k}).
\end{equation}

The next goal is to remove the "non-holomorphic part" of the series expansion of $\Psi(z^{1/k}).$ Define the function $F$ by
\begin{equation}\label{defin of F}
    F(z^k):=\frac{1}{k}\sum_{j=0}^{k-1}\Psi(\zeta^jz),
\end{equation}
where $\zeta:=e^{2\pi i \frac{1}{k}}.$ We shall show that $F$ is a well defined holomorphic function close to $0.$ Indeed, it is enough to check the following finite sums:
\begin{align}
\frac{1}{k}\sum_{j=0}^{k-1}\sum_{q=1}^{Q}\frac{(-1)^{q+1}}{q}\Big(\sum_{l=0}^{L}\lambda_l(\zeta^jz)^{2k+l}\Big)^q,
\end{align}
where $1-\Phi(z)=-\sum_{l=0}^{+\infty}\lambda_lz^{2k+l}.$ Recall the average equality 
\begin{equation}\label{average type}
  \frac{1}{k}\sum_{j=0}^{k-1}\zeta^{jm}=
\begin{cases}
1,  &\text{ if } k \text{ $\mathbb Z$-divides }m\\
0, &\text{ otherwise. }
\end{cases}  
\end{equation}

The multinomial formula and a straightforward computation yield
\begin{align*}
\frac{1}{k}\sum_{j=0}^{k-1}&\Big(\sum_{l=0}^{L}\lambda_l(\zeta^jz)^{2k+l}\Big)^q=\frac{1}{k}\sum_{j=0}^{k-1}\sum_{|\iota|=q}\frac{q!}{\iota!}\lambda_0^{\iota_0}\cdots\lambda_{L}^{\iota_L}(\zeta^jz)^{2k\iota_0}\cdots(\zeta^jz)^{(2k+L)\iota_L}\\
&=\sum_{|\iota|=q}\frac{q!}{\iota!}\lambda_0^{\iota_0}\cdots\lambda_{L}^{\iota_L}\Big[\frac{1}{k}\sum_{j=0}^{k-1}\zeta^{j(2k\iota_0+(2k+1)\iota_1+...+(2k+L)\iota_L)}\Big]z^{2k\iota_0+(2k+1)\iota_1+...+(2k+L)\iota_L}\\
&=\sum_{|\iota|=q}\frac{q!}{\iota!}\lambda_0^{\iota_0}\cdots\lambda_{L}^{\iota_L}\Big[\frac{1}{k}\sum_{j=0}^{k-1}\zeta^{j(2kq+\iota_1+...+L\iota_L)}\Big]z^{2kq+\iota_1+...+L\iota_L}.
\end{align*}
An application of the equality~\eqref{average type} shows that the non-zero coefficients of the monomials $z^{2kq+\iota_1+...+L\iota_L}$ in the computation above, are those that $k$ $\mathbb Z$-divides the positive integer $2kq+\iota_1+...+L\iota_L.$ Summing up, we reach a classification of the holomorphic function $F:$
\begin{equation}\label{class F}
F(z)=\lambda_0z^2+o(z^{2}), \text{ and } F((\zeta^jw)^k)=F(w^k), \quad j=0,1,...,k-1,
\end{equation}
for all $z,w$ close to $0,$ where $o(z^{2})$ is a holomorphic function and its power series expansion consists of terms $c_mz^{m},$ $m\geq 3.$ The second equality in \eqref{class F} follows immediately from the classification of $F$ and the condition $\zeta^{jmk}=1,$ for all $m=1,2,...$ Let us point out that the series expansion of the difference $\Psi(z)-F(z^k)$ consists only of terms $c_lz^{2k+l},$ where $l\notin\{0,k,2k,...\};$ averaging $\Psi$ actually cancels the rest of the terms in this difference.

Next, define the real analytic function $u$ and the holomorphic function $\delta,$ close to $0,$ by
\begin{equation}\label{def u and delta}
    u(z):=\sigma(|z|^2)-\textrm{Re }F(z), \text{ and } \delta(z):=\Psi(z)-F(z^k),
\end{equation}
respectively.  We obtain,
\begin{align*}
u(z)&=\frac{1}{k}\sum_{j=0}^{k-1}\sigma(|\zeta^{j}w|^{2k})- \frac{1}{k}\sum_{j=0}^{k-1}\textrm{ Re }\Psi(\zeta^{j}w)
=\frac{1}{k}\sum_{j=0}^{k-1}[\sigma(|\zeta^{j}w|^{2k})-\textrm{ Re }\Psi(\zeta^{j}w)]\geq 0,
\end{align*}
where $w^k=z,$ $j=0,1,...,k-1,$ and $z$ close to $0$; the non-negativity comes from \eqref{real psi ineq sigma}. Moreover, applying \eqref{class F} and \eqref{real psi ineq sigma}, we see that
\begin{equation}\label{u ineq re delta}
\begin{aligned}
u(z)&=\sigma(|\zeta^{j}w|^{2k})-\textrm{ Re }\Psi(\zeta^jw)+\textrm{ Re }\Psi(\zeta^jw)-\textrm{Re }F(w^k)\\
&\geq\textrm{ Re }[\Psi(\zeta^jw)-F((\zeta^jw)^k)]\\
&=\textrm{ Re }\delta(\zeta^jw).
\end{aligned}
\end{equation}

We shall classify $u.$ Note that $\partial_{z\bar{z}}u(0)>0,$  hence, one can therefore rotate coordinates in $\CC$ so that
\begin{equation*}
    u(z)=ax^2+by^2+o(|z|^2),
\end{equation*}
where $z:=x+iy$ and $b>0,$ see \cite{HW}. Whence, $\partial_{yy}u(0)=2b\neq 0,$ and the real analytic version of implicit function theorem applied to $\partial_yu$, yields a unique real analytic function $b_0(x)$ with $b_0(0)=b'_0(0)=0$, see \cite{Krantz}. In particular,
$b_0(x)=O(x^2).$ Set
$a_0(x):=u(x,b_0(x))\geq 0.$ By the Taylor expansion with respect to $\alpha_0$ and $b_0,$ we obtain
\begin{equation}\label{classifi u}
    u(z)=a_0(x)+(y-b_0(x))^2\epsilon(z),
\end{equation}
for all $z$ close enough to $0,$ where $\epsilon$ is real analytic, and $\epsilon(0)=b.$

We shall show that $\delta\equiv 0$ utilizing inequality \eqref{u ineq re delta}. This identification leads to $\Psi(z)\equiv\Phi(z^k),$ thereby proving that $h$ is actually an holomorphic function.

Indeed, let us suppose in contrary that $\delta$ is not identically zero. The Taylor power series expansions of $\Psi$ and $F,$ together with the assumptions $k>1$ and $\lambda_l\neq 0,$ for at least one $l\notin\{ 0,k,2k,3k,...\},$  yield
\begin{equation}\label{classif delta}
\delta(z)=\lambda_*z^{\ell}+o(z^\ell),
\end{equation}
where $\lambda_*\in \{\lambda_l:l\notin\{0,k,2k,3k,...\}\},$ $\lambda_*\neq 0,$ and the positive integer $\ell$ cannot be $\mathbb Z$-divided by $k$ and it satisfies $\ell\geq 2k+1.$ Set $r:=\ell/k\in \mathbb Q,$ and 
\begin{equation*}
M:=\max_{j=0,...,k-1}\textrm{Re}(\lambda_*e^{2\pi ijr}), \text{ and } \tilde{M}:=\max_{j=0,...,k-1}\textrm{Re}(\lambda_*e^{i\pi r}e^{2\pi ijr}).
\end{equation*}
Then $M>0$ or $\tilde{M}>0.$ Indeed, the average equality~\eqref{average type} applied for $e^{2\pi ijr}=\zeta^{j\ell},$ shows that 
\begin{equation}\label{M}
   \frac{1}{k}\sum_{j=0}^{k-1}\textrm{Re}(\lambda_*e^{2\pi ijr})=\frac{1}{k}\sum_{j=0}^{k-1}\textrm{Re}(\lambda_*e^{i\pi r}e^{2\pi ijr})=0; 
\end{equation} this forces $M$ to be non-negative. Assume that $M=0;$ we shall show that $\tilde{M}>0.$ If $\textrm{Re}(\lambda_*e^{2\pi ijr})\neq 0,$ for some $j,$ then, by \eqref{M},  there exists at least one $j'$ (possibly $j'=j$) such that $\textrm{Re}(\lambda_*e^{2\pi ij'r})>0,$ which contradicts the assumption $M=0.$ Hence, $M=0$ yields $\textrm{Re}(\lambda_*e^{2\pi ijr})=0$ for all $j=0,...,k-1.$ In particular, for $j=0$: $\lambda_*=\textrm{Im }\lambda_*\neq 0,$ and for $j=1$: $\textrm{Re }(\lambda_*e^{2\pi ir})=\textrm{Im }\lambda_*\sin(2\pi r)=0.$
Whence, $r=q/2,$ for some $q=1,3,5,...,$ and $e^{i\pi r}=\pm i.$ Ergo, equation \eqref{M} yields $\tilde{M}>0.$

As we have already mentioned, we shall utilize \eqref{u ineq re delta} to get contradiction, so we should study the behaviour of $u$ and $\delta$ near $0\in \CC.$ 

First, we deal with $\textrm{Re }\delta(e^{2\pi i\frac{j}{k}}z^{\frac{1}{k}}),$ as $z\rightarrow 0$ through a proper path. This proper path comes from the application of  implicit function theorem. Recall that $b_0(x)=O(x^2),$ and put $z=z(x):=x+ib_0(x)=x(1+O(x)).$ Letting $x>0,$ we may substitute $\zeta^jz^{1/k}=(e^{2\pi ij}z)^{1/k},$ $j=0,1,...,k-1,$ to the Taylor expansion \eqref{classif delta} of the holomorphic function  $\delta$ (picking first a proper branch-cut around zero; this is possible since $(x,b_0(x)),$ $x\geq0,$ yields a real analytic curve on the half complex plain reaching the end point $0\in \CC$), namely,
\begin{equation}\label{delta x pos}
\delta(\zeta^jz^{1/k})=x^r[\lambda_*e^{2\pi ijr}(1+O(x))^r+o(x^s)],
\end{equation}
for all $x>0$ close to $0,$ and some $s>0.$ Similarly, if $x<0,$ then $x=|x|e^{i\pi},$ and thus,
\begin{equation}\label{delta x neg}
\delta(\zeta^jz^{1/k})=|x|^r[\lambda_*e^{i\pi r}e^{2\pi ijr}(1+O(|x|))^r+o(|x|^s)],
\end{equation}
for all $x<0$ close to $0.$ Recall that $M>0$ or $\tilde{M}>0.$ Without loss of generality, let us suppose that $M>0$ (we may follow the same procedure for $\tilde{M}$ since $|x|>0.$) Then there exists $j\in \{0,...,k-1\}$ such that $\textrm{Re}(\lambda_*e^{2\pi ijr})>0.$ For this particular $j,$ we have then
$$\textrm{Re}\Big(\lambda_*e^{2\pi ijr}(1+O(x))^r+o(x^s)\Big)\longrightarrow \textrm{Re}(\lambda_*e^{2\pi ijr})>0, \quad \text{ as } x\rightarrow 0^+.$$
Equations \eqref{delta x pos} and \eqref{delta x neg}, together with the fact $M>0$ or $\tilde{M}>0,$ yield
\begin{equation}\label{estimate re delta xr}
\textrm{Re }\delta(\zeta^jz^{1/k})\geq c|x|^r,  
\end{equation}
for some fixed $j,$ and for all $x$ close to $0,$ either by the side $x\geq 0,$ or the side $x\leq 0.$ The constant $c>0$ is uniform with respect to $x.$ 

In what follows, each constants $c,C>0$ are uniform with respect to the variables $z,w,x,y$. 

Summing up, the inequalities \eqref{u ineq re delta}, \eqref{estimate re delta xr} hold for all $z=x+ib_0(x),$ $w^k=z,$ where either $x\geq 0,$ or $x\leq 0,$ close enough to $0.$ For simplicity choose the side $x\geq 0.$ Since $a_0(x)=u(x,b_0(x))\geq 0,$ we have
\begin{equation}\label{a0 ineq x r}
a_0(x)\geq c|x|^r,
\end{equation}
for all $x\geq 0$ close to $0.$ This shows that the real analytic function $a_0$ is not identically zero. Let $a_0(x)=d_mx^m+d_{m+1}x^{m+1}+...,$ where $d_m\neq 0,$ and $m=1,2,....$ Then 
\begin{equation}\label{a0 teliko}
    a_0(x)\geq c|x|^m,
\end{equation}
for all $x$ close to $0$ (that is, in both sides $x\geq 0$ and $x\leq 0$). More specifically, $m$ must be even due to the positivity of $a_0,$ $d_m> 0,$ and the inequality \eqref{a0 ineq x r} together with the condition $r\in \mathbb Q,$ yield $m<r.$

Furthermore, setting $\beta(x,y):=y-b_0(x),$ and recalling that $\epsilon(0)=b>0,$ we obtain, by \eqref{a0 teliko} and \eqref{classifi u}, that
\begin{equation}\label{u teliko}
    u(z)=a_0(x)+\beta(x,y)^2\epsilon(z)\geq c(|x|^m+\beta(x,y)^2),
\end{equation}
for all $z=x+iy$ close to $0.$

On the other hand, \eqref{classif delta} yields $|\delta(w)|\leq c|z|^r,$ for all $z=w^k.$ In particular, since $y=b_0(x)+\beta(x,y),$ $b_0(x)=O(x^2),$ and thus  $$z=x+iy=x+ib_0(x)+i\beta(x,y)=x(1+O(x))+i\beta(x,y),$$
we get
\begin{equation}\label{delta telika}
|\delta(w)|\leq C(|x|^r+|\beta(x,y)|^r),
\end{equation}
for all $w$ close enough to $0.$ Recall that $r>m$ and $r>2.$

Summarizing, the definitions \eqref{def u and delta} yield the equality
\begin{equation}\label{last iqe}
\begin{aligned}
u(w^k)-\textrm{Re } \delta(w)&=\sigma(|w|^{2k})-\textrm{Re } F(w^k)-\textrm{Re }\Psi(w)+\textrm{Re }F(w^k)\\
&=\sigma(|w|^{2k})-\textrm{Re }\Psi(w),
\end{aligned} 
\end{equation}
for all $w^k=z$ close to $0,$ and the inequalities \eqref{u teliko}, \eqref{delta telika} yield
\begin{equation*}
u(w^k)-\textrm{Re } \delta(w)\geq u(w^k)-|\delta(w)|\geq c(|x|^m+\beta(x,y)^2)-C(|x|^r+|\beta(x,y)|^r),
\end{equation*}
for all $z=x+iy$ and $z=w^k,$ where $x$ and $y$ close enough to $0\in \RR.$ Since $r>m$ and $r>2,$ we obtain that 
\begin{equation}\label{lll}
u(w^k)-\textrm{Re } \delta(w)>0,
\end{equation}
for all non-zero $z$ and $z=w^k$ sufficiently close to $0.$

Last, the proofs of \eqref{h Nash curve}, \eqref{real psi ineq sigma}, \eqref{last iqe}, and \eqref{lll} lead to contradiction on the assumption that $\delta$ is not identically $0.$

It turns out that $\delta\equiv0,$ and hence, $\Psi(z)=F(z^k),$ for all $z$ in a fixed neighbourhood of $0.$ In particular, 
$$h(z)=e^{F(z)},$$
namely, there exists $\epsilon>0$ such that $h\in \Hol(\Delta(0,\epsilon)).$ The assertion follows.

\end{proof}
\end{lemma}

In light of the Lemma~\ref{Lemma}, we achieve a characterization of the boundary zero set of a stable polynomial in $\BB$ (which is a stronger version of \cite[Theorem 1]{VZ}). See also \cite[Proposition 1.3-Theorem 1.5]{Knese} for a different approach.

\begin{theorem}\label{theorem}
Let $p\in \CC[z,w]$ be a stable polynomial such that $\mathcal{Z}(p)\cap \SSS\neq \emptyset.$ Then $\mathcal{Z}(p)\cap \SSS$ is either a finite set or, the non-isolated part of $\mathcal{Z}(p)\cap \SSS$ coincide with $\cup_iM_i,$ where $M_i,$ $i=1,2,..,k,$ $k\geq 2,$ are Nash submanifold lying in $\SSS$, each equivalent to $(-1,1)$ by a complex tangential Nash diffeomorphism. Each $M_i$ contains exactly two points on its boundary $\overline{M_i}\setminus M_i.$ In particular, for each $i$ and each $\alpha\in \overline{M_i}\setminus M_i,$ there exists $j\neq i$ such that $\alpha\in M_j,$ $\alpha\in \overline{M_i\cap M_j},$ and $M_i\cap M_j$ is connected open set in the relative topology of $M_i.$ Furthermore, each Nash diffeomorphism $\phi_i:(-1,1)\rightarrow M_i$--under the application of a unitary matrix $\mathcal{U}_i$ with respect to $\alpha=(0,1)$--is locally, around $0\in (-1,1),$ classified as follows:
$$\mathcal{U}_i\phi_i(t)=(\beta_i t+o(t),1-\frac{|\beta_i|^2}{2}t^2+o(t^2)),$$
for all $t\in (-\epsilon,\epsilon),$ where $\beta_i\in \CC\setminus\{0\},$ and $\epsilon>0$ fixed.
\begin{proof}
By the discussion in Section~\ref{Prel Remarks}, if $\dim_\RR\mathcal{Z}(p)\cap \SSS=1,$ then the non-isolated part coincides with $\cup_iM_i,$ where $M_i,$ are finitely many, pairwise disjoint, complex tangential Nash submanifold lying in $\SSS$, each Nash diffeomorphic to $(-1,1)^{m_i},$ $m_i=0,1.$ In particular, each one-dimensional component $M_i$ has exactly two points on its boundary $\overline{M_i}\setminus M_i.$ Since the components are homeomorphic to $(-1,1),$ there exists at least one $M_{i_0}$ with $m_{i_0}=0.$ This case $m_{i_0}=0$ yields just an accumulation point, namely, the compact manifold $M_{i_0}=\{\alpha\}$. In particular, $\alpha\in \overline{M_{j_0}}\setminus M_{j_0},$ for some fixed $j_0\neq i_0.$ Applying the proof of Lemma~\ref{Lemma} to this accumulation point $\alpha$ and the arguments concerning the complex tangential property in the proof of \cite[Theorem 1]{VZ}, we obtain a Nash submanifold which is equivalent to $(-1,1)$ by a complex tangential Nash diffeomorphism, lies in $\mathcal{Z}(p)\cap\SSS,$ $\alpha\in N_\alpha,$ $\alpha\in \overline{N_\alpha\cap M_{j_0}},$ and $N_\alpha\cap M_{j_0}$ is connected open set in the relative topology of $M_{j_0}.$  Last, for each $M_i$ with $m_i=1,$ we have that  $\cup_{i,m_i=1}(\overline{M_i}\setminus M_i)$ is a finite subset of $\mathcal{Z}(p)\cap\SSS.$ Each point $\alpha\in \cup_{i,m_i=1}(\overline{M_i}\setminus M_i)$ lies on a finite number of $\overline{M_i}\setminus M_i$ (probably more than two), and hence, applying the previous argument to each one of them separately, it lies on a finite number of such $N_\alpha.$ Add the finitely many $N_\alpha$ on the partition $\{M_i\}_{i,m_i=1}.$ Last, the classification of the Nash diffeomorphisms follows from the complex-tangential property, see \cite{VWZ}.
\end{proof}
\end{theorem}

\section{Further discussion}

\begin{remark}
We note that \cite[Theorem 1.5~(1)]{Knese} singles out the curve-type boundary zero case, establishing that the underlying (complex) algebraic curve intersects the boundary $\SSS$ along a real one-dimensional analytic curve. Hence, Lemma~\ref{Lemma} serves as an alternative, independent proof of the existence of such a curve lying in the intersection of a complex algebraic curve (associated with a stable polynomial) with the boundary sphere. In particular, Lemma~\ref{Lemma} explicitly guarantees that this real analytic curve is non-singular (diffeomorphism onto its image), providing additional geometric regularity to the boundary trace.
\end{remark}

\begin{remark}
Roughly speaking Theorem~\ref{theorem} says that an one-dimensional boundary zero set of a stable polynomial is either a union of finitely many isolated points and finitely many one-dimensional compact manifolds in $\SSS$ or, a union of finitely many isolated points and finitely many one-dimensional submanifolds in $\SSS$ that "cross each other" (i.e., "cusps" cannot occur; see Section~\ref{Prel Remarks} about the term "cusp").  By the term "cross each other", we mean that locally--around a singular point $\alpha\in \mathcal{Z}(p)\cap \SSS$--the boundary zero set mirrors the novel example $i-1+2z_2-(i+1)z_2^2+z_1^2$ of \cite{Knese} around its singular point; with probably an arbitrary large number of finitely many half-branches and not just four. This example serves as a crucial stable polynomial; beyond the well known polynomial $1-z^2-w^2$--and its variations by a unitary matrix--whose boundary zero set is just a compact manifold in $\SSS.$ The existence of such polynomials is essential and goes back to the discussion in \cite[Chapter 10]{Rudin} and the references within. 
\end{remark}

\begin{remark}
Denote the open ball with center $\alpha\in \CC^2$ and radial $\epsilon>0$ by $\Ba(\alpha,\epsilon):=\{(z,w)\in \CC^2:||(z,w)-\alpha||<\epsilon\}.$ Note that Lemma~\ref{Lemma} together with Theorem~\ref{theorem} rule out the existence of a hypothetical closed (real) algebraic curve $\gamma_\alpha([-1,1])\subset \mathcal{Z}(p)\cap \SSS$ and $\epsilon_0>0,$ such that $\gamma_\alpha(0)=\alpha,$ and for each $0<\epsilon<\epsilon_0,$ the curve satisfies $\gamma_\alpha([-1,1])\cap \Ba(\alpha,\epsilon)=(M_0\cup M_1\cup M_2\cup M_3\cup M_4)\cap \Ba(\alpha,\epsilon),$ where $M_0=\{\alpha\},$ $M_i$ are pairwise disjoint submanifolds with $\dim_\RR M_i=1,$ $i=1,2,3,4,$ and $\cap_i\overline{M_i}=M_0,$ and $M_0\cup M_1\cup M_2$ yields a submanifold while $M_0\cup M_3\cup M_4$ does not (i.e. $\overline{M_3}\cup\overline{M_4}$ builds a "cusp" at $\alpha$). Although the point $\alpha$ lies on the submanifold $M_0\cup M_1\cup M_2$ which is contained in $\mathcal{Z}(p)\cap \SSS,$ it also lies on $M_0\cup M_3\cup M_4$ which is contained in $\mathcal{Z}(p)\cap \SSS.$ Hence, Lemma~\ref{Lemma} together with Theorem~\ref{theorem} serve as a stronger result than proving that for each $\alpha\in\mathcal{Z}(p)\cap \SSS$ there exists a Nash curve $\gamma_\alpha:(-1,1)\rightarrow \mathcal{Z}(p)\cap \SSS,$ that is diffeomorphism onto it image and satisfies $\gamma_\alpha(0)=\alpha.$ See Remark 9.5.8 together with Subsections 9.3 and 9.5 of \cite{Algebraic} on the half-branches of algebraic curves. For a local description, Theorem~\ref{theorem} says that for each point $\alpha\in \mathcal{Z}(p)\cap\SSS,$ there exist finitely many one dimensional Nash submanifolds $N_1,...,N_k\subset \mathcal{Z}(p)\cap \SSS$ and $\epsilon_0>0$ such that $\alpha\in N_i,$ for all $i=1,...,k,$ and $\mathcal{Z}(p)\cap\SSS\cap \Ba(\alpha,\epsilon)=(N_1\cup\cdots \cup N_k)\cap \Ba(\alpha,\epsilon),$ for all $0<\epsilon<\epsilon_0.$ 
\end{remark}

\begin{remark}
Recalling the \cite[Definition 4.3]{VWZ} of essential singularities and weak essential singularities (these definitions coincide in the setting $\BB$), Theorem~\ref{theorem} proves that given a stable polynomial $p\in \CC[z,w],$ the set $\mathcal{Z}(p)\cap \SSS$ has no essential singularities.  As we have already mentioned in the introduction, Lemma~\ref{Lemma} and Theorem~\ref{theorem} cannot be extended to the setting $\mathbb B_n,$ for $n\geq 3,$ see \cite{VWZ}. Recall also that a compact set $K\subset \SSS$ is called a \emph{peak set} for $A^\infty(\BB)$ (the space of all holomorphic functions in $\BB$ that are extended smoothly on $\overline{\BB}$), if there exists a  function $f\in A^\infty(\BB)$ such that $f(z)=0$ for all $z\in K,$ and $\textrm{Re } f(z)>0,$ for all $z\in \overline{\BB}\setminus K,$ see \cite[Chapter 10]{Rudin} for more information. The theory established in \cite{VWZ} together with Theorem~\ref{theorem} shows that $\mathcal{Z}(p)\cap \SSS=\cup_iM_i,$ where each $\overline{M_i}$ is a peak set for $A^\infty(\BB).$ Thus, Theorem~\ref{theorem} together with \cite[Proposition 10.7.5-Theorem 10.7.6]{Rudin} yield, additionally, an equivalent characterization of $\mathcal{Z}(p)\cap \SSS$ in terms of peak sets for $A^\infty(\BB).$ Topics concerning peak sets, the boundary zero set of stable polynomials, and characterization of cyclicity of multipliers, are also discussed in \cite{Ziarati}. The stability condition in the proof of Lemma~\ref{Lemma} may be formulated in the switched setting coming up in \cite{Knese}. Theorem~\ref{theorem} together with established theory in \cite{Ziarati} can be applied in the cyclitiy problem of polynomials in Dirichlet-type spaces in $\BB,$ see also \cite{KV}. Lemma~\ref{Lemma} should be also compared with \cite[1.5 Lemma]{Noel}. 
\end{remark}

\begin{remark}
The involvement of Nash curve selection lemma in proof of Lemma~\ref{Lemma} can be omitted; the hypotheses that $\alpha\in \overline{M_i}\setminus M_i$ is enough. We use it to stress the distinction.
\end{remark}

\begin{remark}
One may consider the same problem for polynomials that are stable with respect to domains whose topological boundary is a an algebraic set or a real analytic set. Recall the semi-analytic curve selection lemma in \cite[Lemma 2.28]{Denk} as it is stated by Z. Denkowska and M.P. Denkowski. A simple example is the "ellipsoides" $\Omega:=\{(z,w)\in \CC^2:\sigma(z,w)<1\},$ where $\sigma(z,w):=|z|^{2p}+|w|^{2p},$ where $p>1$ is integer. 
\end{remark}

\section*{Acknowledgements}
The author would like to thank M. Denkowski, J.F. Fernando Galván, \L. Kosi\'{n}ski,  W. Kucharz for enlightenments. The conceptualisation of the main ideas took place during the Bench Math Session workshop, Kraków, February 2025 together with K. Zarvalis. The author would like also to thank D. Bianchi for help with verifying the arguments.

\section*{Statements and Declarations}
\subsection*{Funding}
Not applicable.

\subsection*{Competing Interests}
The author has no competing interest to declare that is relevant to the content of this article.

\subsection*{Data Availability Statement}
Not applicable.

\subsection*{AI Declaration}
ChatGPT 6 assisted for verifying the arguments concerning the classification of the functions $h_j$ in the proof of Lemma~\ref{Lemma}. The text, the ideas, and the strategies are based on human and the references contained herein. The author take full responsibility for the content of this manuscript.

\bibliographystyle{plain}

\end{document}